\documentclass[pdflatex,sn-mathphys-num]{sn-jnl}

\usepackage{amsmath,amssymb,amsthm}

\newtheorem{theorem}{Theorem}
\newtheorem{lemma}{Lemma}
\newtheorem{proposition}{Proposition}
\newtheorem{corollary}{Corollary}

\begin{document}

\title{A Korselt Characterization of Carmichael Numbers of the Form
$(aC_n+1)(aC_{n+1}+1)(aC_{n+2}+1)$}

\author{\fnm{Hayat} \sur{Bensella}}\email{hbensella@usthb.dz}

\affil{\orgdiv{Faculty of Mathematics, LATN Laboratory}, \orgname{USTHB}, \orgaddress{\street{BP 32, El Alia}, \city{Bab Azzouar}, \postcode{16111}, \state{Algiers}, \country{Algeria}}}

\abstract{
We study products
\[
N=(aC_n+1)(aC_{n+1}+1)(aC_{n+2}+1),
\]
where $(C_n)$ is the Lucas--balancing sequence. Under the assumption
that all three factors are prime, we give a necessary and sufficient
condition for $N$ to be a Carmichael number: this happens precisely
when $n$ is even and $a=kC_{n+1}$ for a positive integer $k$ satisfying
\[
48k\equiv 7\pmod{C_nC_{n+2}}.
\]
We derive further congruence restrictions, an explicit parametrization
of the integer solutions of the resulting congruence, a lower bound
for $k$, and a divisibility relation involving the middle prime factor.
The characterization is conditional on the primality of the three
factors and does not assert that such prime triples exist for infinitely
many parameters.
}



\maketitle

\section{Introduction}

A composite positive integer $N$ is called a Carmichael number if
\[
b^{N-1}\equiv 1\pmod N
\]
for every integer $b$ relatively prime to $N$. We use Korselt's
criterion: a composite positive integer $N$ is a Carmichael number
if and only if it is square-free and
$p-1\mid N-1$ for every prime divisor $p$ of $N$.

The Lucas--balancing numbers were introduced by Panda \cite{Panda}
in terms of the balancing numbers $(B_n)_{n\geq 0}$, originally
studied by Behera and Panda \cite{BeheraPanda}. They are defined by
$C_n=\sqrt{8B_n^2+1}$ and satisfy the recurrence
\[
C_0=1,\qquad C_1=3,\qquad
C_{n+2}=6C_{n+1}-C_n,\quad n\geq 0.
\]
Carmichael numbers have been studied since the early works of
Korselt \cite{Korselt1899} and Carmichael \cite{Carmichael1910}.
Erd\H{o}s and Odlyzko \cite{Erdos1979} investigated the density of
odd integers related to primes and powers of $2$. A major
breakthrough came with Alford, Granville, and Pomerance
\cite{Alford1994}, who proved that infinitely many Carmichael
numbers exist. Subsequent studies have addressed their distribution
and arithmetic properties. Harman \cite{Harman2005} investigated
their counting function, while Wright \cite{Wright2012} established
nonexistence results for certain classes of Carmichael numbers.
Further contributions include the work of Banks et al.
\cite{Banks2015} and Cilleruelo, Luca, and Pizarro
\cite{Cilleruelo2016}. More recently, Alahmadi and Luca
\cite{AlahmadiLuca2022} proved that no Carmichael number can be
of the form $2^n p+1$, where $p$ is prime. Advances in computational
enumeration were reported by Shallue and Webster
\cite{ShallueWebster2025}.

In this paper, we apply Korselt's criterion to products of three
consecutive terms of the Lucas--balancing sequence. Throughout, we consider
\[
N=(aC_n+1)(aC_{n+1}+1)(aC_{n+2}+1),
\qquad a,n\geq 1,
\]
assuming that all three factors are prime. Our main result gives an
exact characterization of when $N$ is a Carmichael number.

\begin{theorem}\label{thm:main}
Let $a,n\geq 1$, and suppose that
\[
aC_n+1,\qquad aC_{n+1}+1,\qquad aC_{n+2}+1
\]
are prime. Then
\[
N=(aC_n+1)(aC_{n+1}+1)(aC_{n+2}+1)
\]
is a Carmichael number if and only if $n$ is even and there exists
an integer $k\geq 1$ such that
\[
a=kC_{n+1},
\qquad
48k\equiv 7\pmod{C_nC_{n+2}}.
\]
\end{theorem}

In particular, odd indices are impossible under the three-prime
assumption.

\section{Elementary identities}

We establish several elementary properties of the sequence.

\begin{lemma}\label{lem:positive}
The sequence $(C_j)_{j\geq 0}$ consists of positive odd integers
and is strictly increasing.
\end{lemma}

\begin{proof}
The initial terms satisfy $C_1>C_0>0$. If
$C_{j+1}>C_j>0$, then
\[
C_{j+2}=6C_{j+1}-C_j
>5C_{j+1}>C_{j+1}.
\]
Positivity and strict increase follow by induction. Both initial
terms are odd, and the recurrence shows that every subsequent term
is odd.
\end{proof}

\begin{lemma}\label{lem:identities}
For every $j\geq 0$,
\[
\gcd(C_j,C_{j+1})=1
\]
and
\[
C_jC_{j+2}-C_{j+1}^2=8.
\]
\end{lemma}

\begin{proof}
The recurrence gives
\[
\gcd(C_{j+1},C_{j+2})
=\gcd(C_{j+1},6C_{j+1}-C_j)
=\gcd(C_{j+1},C_j).
\]
Since $\gcd(C_0,C_1)=1$, the first assertion follows.

For the second assertion, put
\[
E_j=C_jC_{j+2}-C_{j+1}^2.
\]
Using $C_{j+3}=6C_{j+2}-C_{j+1}$, we obtain
\[
\begin{aligned}
E_{j+1}
&=C_{j+1}C_{j+3}-C_{j+2}^2\\
&=C_{j+1}(6C_{j+2}-C_{j+1})-C_{j+2}^2\\
&=(6C_{j+1}-C_{j+2})C_{j+2}-C_{j+1}^2\\
&=C_jC_{j+2}-C_{j+1}^2\\
&=E_j.
\end{aligned}
\]
Hence $E_j$ is constant. Since $C_2=17$,
\[
E_0=C_0C_2-C_1^2=17-9=8.
\]
\end{proof}

\begin{lemma}\label{lem:mod3}
For every $j\geq 0$,
\[
3\mid C_j\quad\Longleftrightarrow\quad j\text{ is odd}.
\]
Moreover,
\[
\gcd(C_j,C_{j+2})=
\begin{cases}
1,&j\text{ is even},\\
3,&j\text{ is odd}.
\end{cases}
\]
\end{lemma}

\begin{proof}
Modulo $3$, the recurrence becomes
\[
C_{j+2}\equiv -C_j\pmod 3.
\]
Starting from $C_0\equiv 1$ and $C_1\equiv 0\pmod 3$, we obtain
\[
C_{2r}\equiv (-1)^r\pmod 3,\qquad
C_{2r+1}\equiv 0\pmod 3.
\]
By Lemma~\ref{lem:identities},
\[
\begin{aligned}
\gcd(C_j,C_{j+2})
&=\gcd(C_j,6C_{j+1})\\
&=\gcd(C_j,6).
\end{aligned}
\]
Since every $C_j$ is odd, the claimed formula follows.
\end{proof}

\begin{lemma}\label{lem:mod8}
For every $j\geq 0$,
\[
C_{2j}\equiv1\pmod8,
\qquad
C_{2j+1}\equiv3\pmod8.
\]
\end{lemma}

\begin{proof}
The assertions hold for $j=0$, since $C_0=1$ and $C_1=3$.
Assume that
\[
C_{2j}\equiv1\pmod8,
\qquad
C_{2j+1}\equiv3\pmod8.
\]
Then, using the recurrence,
\[
C_{2j+2}=6C_{2j+1}-C_{2j}
\equiv18-1\equiv1\pmod8,
\]
and
\[
C_{2j+3}=6C_{2j+2}-C_{2j+1}
\equiv6-3\equiv3\pmod8.
\]
The result follows by induction.
\end{proof}

\section{Proof of the main theorem}

We apply Korselt's criterion. Since
\[
C_n<C_{n+1}<C_{n+2}
\]
by Lemma~\ref{lem:positive}, and since $a\ge1$, the three assumed
prime factors of $N$ are distinct. Hence $N$ is square-free.

Since the three factors of $N$ are distinct primes, $N$ is a
Carmichael number if and only if
\[
aC_n\mid N-1,\qquad
aC_{n+1}\mid N-1,\qquad
aC_{n+2}\mid N-1.
\]
Expanding and collecting terms according to powers of $a$, we obtain
\[
\begin{aligned}
N
={}&(aC_n+1)(aC_{n+1}+1)(aC_{n+2}+1)\\
={}&1+a(C_n+C_{n+1}+C_{n+2})\\
&+a^2(C_nC_{n+1}+C_nC_{n+2}
+C_{n+1}C_{n+2})\\
&+a^3C_nC_{n+1}C_{n+2}.
\end{aligned}
\]
Hence,
\[
\begin{aligned}
N-1=a\bigl[&
C_n+C_{n+1}+C_{n+2}\\
&+a(C_nC_{n+1}+C_nC_{n+2}
+C_{n+1}C_{n+2})\\
&+a^2C_nC_{n+1}C_{n+2}
\bigr].
\end{aligned}
\]

For convenience, define
\[
\begin{aligned}
T={}&C_n+C_{n+1}+C_{n+2}\\
&+a(C_nC_{n+1}+C_nC_{n+2}
+C_{n+1}C_{n+2})\\
&+a^2C_nC_{n+1}C_{n+2}.
\end{aligned}
\]
Thus
\[
\boxed{N-1=aT.}
\]

Since $N-1=aT$ and $a\ge1$, Korselt's criterion is equivalent to
\[
C_n\mid T,\qquad
C_{n+1}\mid T,\qquad
C_{n+2}\mid T. \tag{1}
\]

We now derive these divisibility conditions explicitly.

\medskip
\noindent\textbf{Step 1: The condition involving $C_{n+1}$.}

Reducing $T$ modulo $C_{n+1}$ gives
\[
T\equiv C_n+C_{n+2}+aC_nC_{n+2}
\pmod{C_{n+1}}.
\]
By the recurrence relation,
\[
C_{n+2}=6C_{n+1}-C_n,
\]
so that
\[
C_{n+2}\equiv -C_n\pmod{C_{n+1}}.
\]
It follows that
\[
\begin{aligned}
T
&\equiv C_n-C_n-aC_n^2\\
&\equiv -aC_n^2\pmod{C_{n+1}}.
\end{aligned}
\]
The consecutive terms are coprime:
\[
\gcd(C_n,C_{n+1})=1.
\]
Hence
\[
C_{n+1}\mid T
\quad\Longleftrightarrow\quad
C_{n+1}\mid aC_n^2
\quad\Longleftrightarrow\quad
C_{n+1}\mid a.
\]
Therefore, the middle divisibility condition in (1) is equivalent
to the existence of an integer $k\ge1$ such that
\[
\boxed{a=kC_{n+1}.} \tag{2}
\]

\medskip
\noindent\textbf{Step 2: The conditions involving $C_n$ and $C_{n+2}$.}

Reducing $T$ modulo $C_n$, all terms containing $C_n$ vanish,
and we obtain
\[
T\equiv C_{n+1}+C_{n+2}
+aC_{n+1}C_{n+2}\pmod{C_n}.
\]
Since
\[
C_{n+2}\equiv 6C_{n+1}\pmod{C_n},
\]
we get
\[
T\equiv C_{n+1}(7+6aC_{n+1})
\pmod{C_n}.
\]
As $\gcd(C_n,C_{n+1})=1$, this yields
\[
\boxed{
C_n\mid T
\quad\Longleftrightarrow\quad
C_n\mid 7+6aC_{n+1}.
} \tag{3}
\]

Similarly, reducing $T$ modulo $C_{n+2}$ gives
\[
T\equiv C_n+C_{n+1}+aC_nC_{n+1}
\pmod{C_{n+2}}.
\]
Using
\[
C_n\equiv 6C_{n+1}\pmod{C_{n+2}},
\]
we obtain
\[
T\equiv C_{n+1}(7+6aC_{n+1})
\pmod{C_{n+2}}.
\]
Since $\gcd(C_{n+1},C_{n+2})=1$, it follows that
\[
\boxed{
C_{n+2}\mid T
\quad\Longleftrightarrow\quad
C_{n+2}\mid 7+6aC_{n+1}.
} \tag{4}
\]

Combining (3) and (4), the two outer divisibility conditions are
equivalent to
\[
C_n\mid 7+6aC_{n+1},
\qquad
C_{n+2}\mid 7+6aC_{n+1}. \tag{5}
\]

\medskip
\noindent\textbf{Step 3: The parity of $n$.}

The sequence satisfies the identity
\[
\boxed{C_nC_{n+2}-C_{n+1}^2=8.} \tag{6}
\]
Indeed, this identity is Lemma~\ref{lem:identities}.

Substituting $a=kC_{n+1}$ from (2) into (5), we find
\[
C_n\mid 7+6kC_{n+1}^2,
\qquad
C_{n+2}\mid 7+6kC_{n+1}^2.
\]
By (6),
\[
\begin{aligned}
7+6kC_{n+1}^2
&=7+6k(C_nC_{n+2}-8)\\
&=6kC_nC_{n+2}+7-48k.
\end{aligned}
\]
Therefore, (5) is equivalent to
\[
\boxed{
C_n\mid 48k-7,\qquad
C_{n+2}\mid 48k-7.
} \tag{7}
\]

We now show that $n$ must be even. By Lemma~\ref{lem:mod3},
\[
3\mid C_j\quad\Longleftrightarrow\quad j\text{ is odd}.
\]
If $n$ were odd, then both $C_n$ and $C_{n+2}$ would be divisible
by $3$. By (7), this would imply
\[
3\mid 48k-7.
\]
But
\[
48k-7\equiv -1\pmod 3,
\]
which is impossible. Hence
\[
\boxed{n\text{ is even}.} \tag{8}
\]

\medskip
\noindent\textbf{Step 4: Combining the two congruences.}

For even $n$, both $C_n$ and $C_{n+2}$ are odd. By Lemma~\ref{lem:identities},
\[
\gcd(C_n,C_{n+2})\mid 8.
\]
Since their greatest common divisor is odd, we conclude that
\[
\gcd(C_n,C_{n+2})=1.
\]
Consequently, the two divisibilities in (7) are equivalent to
\[
\boxed{48k\equiv 7\pmod{C_nC_{n+2}}.} \tag{9}
\]

This proves the necessity of all the stated conditions.

\medskip
\noindent\textbf{Step 5: Sufficiency.}

Conversely, suppose that $n$ is even and that
\[
a=kC_{n+1},
\qquad
48k\equiv 7\pmod{C_nC_{n+2}}
\]
for some integer $k\ge1$.

We prove that all three divisibilities in (1) hold.

First, since $C_{n+1}\mid a$, the calculation in Step 1 gives
\[
T\equiv -aC_n^2\equiv 0\pmod{C_{n+1}}.
\]
Thus
\[
C_{n+1}\mid T.
\]

Next, the congruence in (9), together with identity (6), gives
\[
\begin{aligned}
7+6aC_{n+1}
&=7+6kC_{n+1}^2\\
&=6kC_nC_{n+2}+7-48k\\
&\equiv 0\pmod{C_nC_{n+2}}.
\end{aligned}
\]
In particular,
\[
C_n\mid 7+6aC_{n+1},
\qquad
C_{n+2}\mid 7+6aC_{n+1}.
\]
By the equivalences established in Step 2, we obtain
\[
C_n\mid T,\qquad C_{n+2}\mid T.
\]

Therefore,
\[
C_n\mid T,\qquad
C_{n+1}\mid T,\qquad
C_{n+2}\mid T.
\]
Since $N-1=aT$, it follows that
\[
aC_n\mid N-1,\qquad
aC_{n+1}\mid N-1,\qquad
aC_{n+2}\mid N-1.
\]
The three factors of $N$ are distinct primes by hypothesis, so $N$
is square-free. Korselt's criterion now implies that $N$ is a
Carmichael number.

The proof is complete.

\section{Further arithmetic consequences}

\begin{lemma}\label{lem:mod48}
For every $j\geq0$,
\[
C_{2j}\equiv
\begin{cases}
1\pmod{48},&j\text{ is even},\\
17\pmod{48},&j\text{ is odd},
\end{cases}
\qquad
C_{2j+1}\equiv3\pmod{48}.
\]
In particular, for every even $n\geq2$,
\[
C_nC_{n+2}\equiv17\pmod{48}.
\]
\end{lemma}

\begin{proof}
We prove the residue pattern by induction on $j$. For $j=0$,
\[
C_0=1\equiv1\pmod{48},
\qquad
C_1=3\equiv3\pmod{48}.
\]
Assume that
\[
C_{2j}\equiv r\pmod{48},
\qquad
C_{2j+1}\equiv3\pmod{48},
\]
where $r\in\{1,17\}$. Then
\[
C_{2j+2}=6C_{2j+1}-C_{2j}
\equiv18-r\pmod{48}.
\]
Thus $r=1$ gives $C_{2j+2}\equiv17\pmod{48}$, while $r=17$ gives
$C_{2j+2}\equiv1\pmod{48}$. Moreover,
\[
C_{2j+3}=6C_{2j+2}-C_{2j+1}
\equiv6(18-r)-3\equiv3\pmod{48},
\]
because $r\in\{1,17\}$ implies $108-6r\equiv6\pmod{48}$.

Hence, by induction,
\[
C_{2j}\equiv1,17,1,17,\ldots\pmod{48},
\qquad
C_{2j+1}\equiv3\pmod{48}.
\]
For even $n$, the residues of $C_n$ and $C_{n+2}$ are $1$ and
$17$ in some order. Their product is therefore congruent to
$17$ modulo $48$.
\end{proof}

\begin{theorem}\label{thm:quotient}
Under the hypotheses of Theorem~\ref{thm:main}, if $N$ is
Carmichael, then
\[
\boxed{\frac{48k-7}{C_nC_{n+2}}\equiv25\pmod{48}.}
\]
\end{theorem}

\begin{proof}
By Theorem~\ref{thm:main}, $n$ is even, $a=kC_{n+1}$ for some
integer $k\ge1$, and
\[
q=\frac{48k-7}{C_nC_{n+2}}
\]
is an integer. Since $k\ge1$, we have $48k-7>0$; also
$C_nC_{n+2}>0$. Hence $q$ is a positive integer.

By Lemma~\ref{lem:mod48},
\[
C_nC_{n+2}\equiv17\pmod{48}.
\]
Since $48k-7\equiv41\pmod{48}$, we get
\[
17q\equiv41\pmod{48}.
\]
As $17^2\equiv1\pmod{48}$, multiplication by $17$ yields
\[
q\equiv17\cdot41\equiv25\pmod{48}.
\]
\end{proof}

\begin{theorem}[Divisibility restriction on the middle prime]
\label{thm:middle-prime}
Under the hypotheses of Theorem~\ref{thm:main}, suppose that
$N=p_0p_1p_2$ is Carmichael. Put
\[
M=C_nC_{n+2},\qquad
q=\frac{48k-7}{M},
\]
where $a=kC_{n+1}$. Then
\[
\boxed{M\mid 6p_1+1}
\]
and
\[
\boxed{\frac{6p_1+1}{M}=6k-q.}
\]
In particular, the quotient on the left is a positive integer.
\end{theorem}

\begin{proof}
By Theorem~\ref{thm:main}, we have
\[
a=kC_{n+1},\qquad M\mid 48k-7.
\]
The identity
\[
C_{n+1}^2=C_nC_{n+2}-8=M-8
\]
gives
\begin{align*}
6p_1+1
&=6(aC_{n+1}+1)+1\\
&=6kC_{n+1}^2+7\\
&=6k(M-8)+7\\
&=6kM-(48k-7).
\end{align*}
Since $M\mid48k-7$, we obtain $M\mid6p_1+1$.
Dividing the displayed identity by $M$ yields
\[
\frac{6p_1+1}{M}
=6k-\frac{48k-7}{M}
=6k-q.
\]
The quotient is positive because $6p_1+1>0$ and $M>0$.
\end{proof}

\begin{corollary}\label{cor:strongbound}
Under the hypotheses of Theorem~\ref{thm:quotient},
\[
k\geq\frac{25C_nC_{n+2}+7}{48}
\]
and
\[
a\geq
\frac{C_{n+1}(25C_nC_{n+2}+7)}{48}.
\]
\end{corollary}

\begin{proof}
By Theorem~\ref{thm:quotient}, the integer
\[
q=\frac{48k-7}{C_nC_{n+2}}
\]
is positive and congruent to $25$ modulo $48$. Hence $q\ge25$.
Thus
\[
48k-7\geq25C_nC_{n+2}.
\]
Rearranging gives the first inequality; the second follows from
$a=kC_{n+1}$.
\end{proof}

\begin{theorem}\label{thm:parametrization}
Fix an even integer $n\geq2$ and put
\[
M=C_nC_{n+2}.
\]
Every positive integer $k$ satisfying
\[
M\mid48k-7
\]
has the form
\[
\boxed{
k=\frac{25M+7}{48}+tM,\qquad t\in\mathbb Z_{\geq0}.
}
\]
Among these values, $k$ is even if and only if
\[
t\equiv\frac{25M+7}{48}\pmod2.
\]
\end{theorem}

\begin{proof}
By Lemma~\ref{lem:mod48}, $M\equiv17\pmod{48}$. Hence
\[
25M+7\equiv25\cdot17+7\equiv0\pmod{48},
\]
so
\[
k_0=\frac{25M+7}{48}
\]
is a positive integer. Moreover,
\[
k_0=\frac{25M+7}{48}<M.
\]

If $M\mid48k-7$, then
\[
48k-7=qM
\]
for some integer $q$. Reducing modulo $48$ shows that
$q\equiv25\pmod{48}$, so $q=25+48t$ for an integer $t$. It follows
that
\[
k=\frac{7+qM}{48}
=\frac{25M+7}{48}+tM
=k_0+tM.
\]
Since $0<k_0<M$, positivity of $k$ is equivalent to $t\ge0$.

Finally, $M$ is odd, so
\[
k=k_0+tM\equiv k_0+t\pmod2.
\]
Thus $k$ is even if and only if $t\equiv k_0\pmod2$, as claimed.
This parametrizes the positive integer solutions of the congruence
only; the primality of the three corresponding factors must still
be checked separately.
\end{proof}

\begin{proposition}[Further congruence restrictions]\label{prop:congruences}
If $N$ is Carmichael under the hypotheses of Theorem~\ref{thm:main},
then its prime factors $p_i=aC_{n+i}+1$ satisfy
\[
p_i\equiv1\pmod6\qquad(i=0,1,2),
\]
and the two outer factors satisfy
\[
p_0\equiv p_2\pmod8.
\]
\end{proposition}

\begin{proof}
Theorem~\ref{thm:main} gives $a=kC_{n+1}$ with $n$ even. By
Lemma~\ref{lem:mod3}, $3\mid C_{n+1}$, so $3\mid a$. Also, $a$
must be even: if it were odd, then each $aC_{n+i}+1$ would be an
even integer greater than $2$, contradicting primality. Thus $6\mid a$
and
\[
p_i=aC_{n+i}+1\equiv1\pmod6.
\]
Since $n$ is even, Lemma~\ref{lem:mod8} gives
$C_n\equiv C_{n+2}\equiv1\pmod8$. Therefore
\[
p_0=aC_n+1\equiv a+1\equiv aC_{n+2}+1=p_2\pmod8.
\]
\end{proof}

\begin{corollary}\label{cor:power2}
Under the three-prime assumption, if $a$ is a positive power of
$2$, then $N$ is not Carmichael.
\end{corollary}

\begin{proof}
Theorem~\ref{thm:main} requires $C_{n+1}\mid a$. But $C_{n+1}$
is odd and greater than $1$ for $n\geq1$, so it cannot divide a
power of $2$.
\end{proof}

\section{An explicit example}\label{sec:example}

Take $n=2$. Then
\[
C_2=17,\qquad C_3=99,\qquad C_4=577,
\]
and
\[
M=C_2C_4=17\cdot577=9809.
\]
The congruence $48k\equiv7\pmod{9809}$ is equivalent to
\[
k\equiv5109\pmod{9809}.
\]
Since $k$ must be even, we may write
\[
k=14918+19618s,\qquad s\geq0.
\]
Taking $s=363$ gives
\[
k=7136252,\qquad
a=99k=706488948.
\]
The corresponding factors are
\[
\begin{aligned}
p_0&=aC_2+1=12010312117,\\
p_1&=aC_3+1=69942405853,\\
p_2&=aC_4+1=407644122997.
\end{aligned}
\]
For this example, the divisibility quotient is
\[
\frac{48k-7}{C_2C_4}
=\frac{342540089}{9809}=34921,
\]
and $34921\equiv25\pmod{48}$, as required by
Theorem~\ref{thm:quotient}.

The primality of the three displayed factors is verified in
Appendix~\ref{app:verification}. This computation concerns only the
particular example and is not used in the proof of
Theorem~\ref{thm:main}.

\appendix

\section{Verification of the example}
\label{app:verification}

For reproducibility, the primality of the three factors in
Section~\ref{sec:example} can be checked using the following
Python/SymPy commands:

\begin{verbatim}
from sympy import isprime

p0 = 12010312117
p1 = 69942405853
p2 = 407644122997

assert isprime(p0)
assert isprime(p1)
assert isprime(p2)
\end{verbatim}

This verifies the primality claims in the displayed example. It does
not form part of the proof of the main theorem.

\backmatter

\section*{Declarations}

\textbf{Conflict of interest.}
The author declares that there is no conflict of interest.

\end{document}